\documentclass[12pt]{article}
\usepackage{amsmath,amssymb}
\usepackage{amsthm,amscd}
\usepackage{graphicx,tikz}
\usepackage{enumerate}
\usepackage{titlesec}
\usepackage{epstopdf}
\usepackage{caption}
\usepackage{xcolor}
\usepackage[sans]{dsfont}
\numberwithin{equation}{section}
\newtheorem{Thm}{Theorem}[section]
\newtheorem*{Thm*}{Theorem}
\newtheorem{Prop}[Thm]{Proposition}
\newtheorem*{Prop*}{Proposition}
\newtheorem{Lem}[Thm]{Lemma}

\newtheorem{Cor}[Thm]{Corollary}
\newtheorem*{Cor*}{Corollary}
\newtheorem{Fact}[Thm]{Fact}

\newtheorem{Exa}[Thm]{Example}

\theoremstyle{remark}

\theoremstyle{definition}

\newtheorem{Def}[Thm]{Definition}

\newtheorem*{Def*}{Definition}

\numberwithin{equation}{section}

\newcommand{\g}[1]{{\mathfrak #1}}
\newcommand{\m}[1]{\mathbb{ #1}}
\newcommand{\mc}[1]{\mathcal{ #1}}

\def\al{\alpha}       \def\be{\beta}        \def\ga{\gamma}
\def\de{\delta}         
       \def\la{\lambda}      
\def\si{\sigma}                
               
\def\om{\omega}       \def\Ga{\Gamma}

\theoremstyle{definition}

\theoremstyle{remark}

\newtheorem{Rmq}[Thm]{Remark}

\numberwithin{equation}{section}
\newfont{\goth}{eufm10 at 12pt}
\newfont{\gots}{eufm8 at 9pt}

\def\bt{\begin{Thm}}
\def\et{\end{Thm}}
\def\br{\begin{Rmq}}
\def\er{\end{Rmq}}

\def\bc{\begin{Cor}}
\def\ec{\end{Cor}}
\def\bp{\begin{Prop}}
\def\ep{\end{Prop}}
\def\bl{\begin{Lem}}
\def\el{\end{Lem}}
\def\bd{\begin{Def}}
\def\ed{\end{Def}}
\def\bq{\begin{quotation}}
\def\eq{\end{quotation}}
\def\bfa{\begin{Fact}}
\def\efa{\end{Fact}}
\def\bexa{\begin{Exa}}
\def\eexa{\end{Exa}}

\def\vs{\vspace{1em}}

\begin{document}
%\mbox{ }\vspace*{-5em}
\title{
On the rational symplectic group
}
\author{Yves Benoist
}
\date{}

\maketitle
%\vspace{-2em}
%\centerline{\footnotesize \mbox{Preliminary version } }

%\hfill { Dedicated to Jacques Tits}
\vs 

\begin{abstract}
\noindent 
%We prove an adapted basis theorem for the rational symplectic group.
We give a short  proof of an elementary classical result:
any rational symplectic matrix can be put in diagonal form 
after right and left multiplication by integral symplectic matrices.
We also give a new proof for its extension to Chevalley groups due to Steinberg by using the Cartan-Bruhat-Tits decomposition over $p$-adic fields.
\end{abstract}

\renewcommand{\thefootnote}{\fnsymbol{footnote}} 
\renewcommand{\thefootnote}{\arabic{footnote}}

%{\footnotesize \tableofcontents}\nopagebreak
%\newpage

%10
\section{Introduction}
\label{secintexp}

In this expository paper 
I present a short proof of a classical theorem 
I needed in \cite{BenoistConvolution2}:
a decomposition of the group ${\rm Sp}(n,\m Q)$
of symplectic matrices with rational coefficients
that gives a parametrization of the double quotient
${\rm Sp}(n,\m Z)\backslash {\rm Sp}(n,\m Q)/{\rm Sp}(n,\m Z)$
where ${\rm Sp}(n,\m Z)$ is 
the subgroup of symplectic matrices with integral coefficients.

This decomposition which can  already be found in \cite{ShimuraArithmetic63} is a symplectic version 
of the ``adapted basis theorem'' 
for $\m Z$-modules, or of the ``Smith normal form'' for integral matrices.
\vs 

In Section \ref{secnorfor} we state precisely this decomposition
that we call the ``symplectic Smith normal form''.

In Section \ref{secanalogy} we explain the analogy with the Cartan-Bruhat-Tits decomposition.

In Section \ref{secbuilding} we recall the relevance of Bruhat-Tits buildings
in this kind of decomposition.

In Section \ref{secconcri} we give an elementary proof of the symplectic 
Smith normal form.

In Section \ref{secstrapp} we give a non-elementary proof of the symplectic Smith normal form that will be applied to other
simply-connected split semisimple algebraic groups {\bf G} defined over $\m Q$ in the last section. 
Indeed we explain how this symplectic Smith normal form 
can be deduced from the Cartan-Bruhat-Tits decomposition 
together with the strong approximation theorem.  

In Section \ref{secchegro} we explain the
extension due to Steinberg of the Smith normal form to the
simply-connected $\m Q$-split groups,
see Theorem \ref{thmhsidsic}.\vs  

The last two sections 
%might be useful for graduate students as 
are a concrete illustration of a classical  strategy: if you want to prove a theorem
over a global field, prove it first over local fields and 
then use a local-global principle.\vs

I would like to  thank Hee Oh for a very helpful comment on a first draft of this note.

%20
\section{The symplectic Smith normal form}
\label{secnorfor}

For any commutative ring $R$ with unity, we denote by  ${\rm Sp}(n,R)$ 
the symplectic group with coefficients in $R$.
This group is the stabilizer of the symplectic form 
$\om $ on $R^{2n}$  given by, for all $x$, $y$ in $R^{2n}$, 
$$
\om(x,y)={}^tx\, J\, y\, 
$$
where 
$
J=\mbox{\scriptsize 
$\left(\!\begin{array}{cc} {\bf 0}&\!\mathds{1}_n\!\\  
\!-\!\mathds{1}_n\!&{\bf 0}\end{array}\!\right)$} .
$ 
Equivalently, one has
$$
{\rm Sp}(n,R):=\{ g\in {\rm GL}(2n,R)\mid \; {}^tg Jg =J\},
$$

If we write the elements of the symplectic group
as  block matrices with blocks of size $n$, one has
$$
{\rm Sp}(n,R)=\{
g=\mbox{\scriptsize 
$\left(\!\begin{array}{cc} \!\al\!&\!\be\!\\  
\!\ga\!&\!\de\!\end{array}\!\right)$} \mid \; 
{}^t\al\ga={}^t\ga\al,\;\;
{}^t\be\de={}^t\de\be,\;\;
{}^t\al\de-{}^t\ga\be=\mathds{1}_n
\}.
$$

\bt
\label{thmhsidsi}
Let $g\in {\rm Sp}(n,\m Q)$.
Then there exist two matrices $\si$ and $\si'$ 
in 
${\rm Sp}(n,\m Z)$ and a positive integral diagonal matrix ${\bf d}={\rm diag}(d_1,\ldots,d_n)$ 
with  $d_1|d_2|\ldots|d_{n}$, and such that
\begin{eqnarray*}
g&=&\si\;
\mbox{\small%\scriptsize 
$\left(\!
\begin{array}{cc} 
{\bf d}&{\bf 0}\\  
{\bf 0}&\!{\bf d}^{-1}\!
\end{array}
\!\right)$}
\;\si'.
\end{eqnarray*}
\et

The condition that the coefficients $d_j$ are positive integers with $d_1$ dividing $d_2$, with
$d_2$ dividing $d_3$, $\ldots$, and $d_{n-1}$ dividing $d_n$ ensures that the diagonal matrix ${\bf d}$ is unique.

I use this precise Theorem \ref{thmhsidsi} as a key tool 
for an apparently completely unrelated problem
in my paper 
%\cite{BenoistConvolution1} and  
\cite{BenoistConvolution2}.
This problem is the construction of 
functions $f$ on the cyclic group $\m Z/d\m Z$ of odd order
whose convolution square is proportional to their square.
Indeed the construction relies on an auxiliary abelian variety
endowed with a unitary $\m Q$-endomorphism $\nu$, 
the symplectic form $\om$ shows up as a polarization of the abelian variety,
and the rational symplectic matrix $g$ 
shows up as the ``holonomy'' of $\nu$.

The first reference to Theorem \ref{thmhsidsi} that  I know is  Shimura's paper \cite[Prop. 1.6]{ShimuraArithmetic63}. 
Moreover in \cite{ShimuraModular63}, Shimura points out the relevance of this theorem to 
show the commutativity of a Hecke algebra and hence to better understand the modular forms on Siegel upper halfspace.
This theorem is also in \cite[p.232]{Freitag83} and is also used by Clozel, Oh and Ullmo in \cite[p.23]{ClozelOhUllmo}.

As we have seen, there is a version of Theorem \ref{thmhsidsi} for the linear group
${\rm SL}(n,\m Q)$, see for instance Proposition \ref{prohsidsi0}.
More generally, there is also a 
version of Theorem \ref{thmhsidsi} for any simply-connected split semisimple algebraic group ${\bf G}$
defined over $\m Q$, if one chooses suitably the $\m Z$-form,
see Section \ref{secchegro}.

%30
\section{The symplectic group over local fields}
\label{secanalogy}
	
Before going on I would like to emphasize the analogy of this theorem with 
two classical theorems. These two classical theorems 
are valid for all
algebraic semisimple groups $G$ and are due respectively to E. Cartan and 
to F. Bruhat and J. Tits. I will not quote here 
their general formulation that can be found respectively in \cite{Helgason01} and in \cite{BruhatTits2}
but only the special case where $G$ is the symplectic group.
\vs 

The first theorem is a decomposition theorem over the real field $\m R$ due to E. Cartan which is called 
either the ``polar decomposition'' or the ``Cartan decomposition''.
We set
\begin{eqnarray*}
{\rm SO}(2n)&:=&\{ g\in {\rm GL}(2n,\m R)\mid \; {}^tg g ={\bf 1}_{2n}\}
\;\;{\rm and}\\
{\rm Sp}(n)&:=&{\rm Sp}(n,\m R)\cap {\rm SO}(2n).
\end{eqnarray*}

Note that the group ${\rm Sp}(n)$ is a maximal compact subgroup of the group
${\rm Sp}(n,\m R)$.

\bt 
{\bf (Cartan)}
\label{thmhsidsir} 
%Let $G:={\rm Sp}(n,\m R)$ and $K={\rm Sp}(n,\m R)\cap {\rm SO}(2n,\m R)$ 
Let $g\in {\rm Sp}(n,\m R)$. Then there exist two matrices $\si$ and $\si'$ in 
${\rm Sp}(n)$ and a positive real diagonal matrix ${\bf d}={\rm diag}(d_1,\ldots,d_n)$ 
with $d_1\leq d_2\leq\ldots \leq d_n\leq 1$  such that
\begin{eqnarray*}
g&=&\si\;
\mbox{\small%\scriptsize 
$\left(\!
\begin{array}{cc} 
{\bf d}&{\bf 0}\\  
{\bf 0}&\!{\bf d}^{-1}\!
\end{array}
\!\right)$}
\;\si'.
\end{eqnarray*}
\et

The second theorem is a decomposition theorem over a local non archi\-medean field $k$ due to F. Bruhat and J. Tits.
We denote by $\mc O_k$ the ring of integers of $k$ and choose a uniformizer $\pi$ in $k$, i.e. a generator of the maximal ideal of $\mc O_k$.
\vs 

Note again that the group ${\rm Sp}(n,\mc O_k)$ is a maximal compact subgroup of the group
${\rm Sp}(n,k)$.

\bt
\label{thmhsidsik} {\bf (Bruhat, Tits)}
%Let $G:={\rm Sp}(n,\m R)$ and $K={\rm Sp}(n,\m R)\cap {\rm SO}(2n,\m R)$ 
%Let $k$ be a non-archimedean local field. 
%$\mc O_k$ be its ring of integers and $\pi\in k$ a uniformizer.
Let $g\in {\rm Sp}(n,k)$. Then there exist two matrices $\si$ and $\si'$ in 
${\rm Sp}(n,\mc O_k)$ and a diagonal matrix ${\bf d}={\rm diag}(\pi^{p_1},\ldots,\pi^{p_n})$ 
with $p_1\geq p_2\geq\ldots\geq p_n\geq 0$ integers  such that
\begin{eqnarray*}
g&=&\si\;
\mbox{\small%\scriptsize 
$\left(\!
\begin{array}{cc} 
{\bf d}&{\bf 0}\\  
{\bf 0}&\!{\bf d}^{-1}\!
\end{array}
\!\right)$}
\;\si'.
\end{eqnarray*}
\et

The analogy between these three theorems is striking. 
It extends the analogy between 
the  Smith normal form of an integral matrix and
the singular value decomposition of a real matrix.
\vs 

In this analogy 
{\it the group of integers points of a group defined over the rational} 
should be handled as
{\it the maximal compact subgroup of a group defined over the real}.
This rough analogy is an equality when dealing with non archimedean
local field.
%at the heart of the point of view of Bruhat and Tits 
%in their construction and their analysis of buildings.
Indeed, when $k$ is a non-archimedean local field,
the group of integer points is an open compact subgroup.

%40
\section{Bruhat-Tits buildings}
\label{secbuilding}

%This paper is dedicated to Jacques Tits who, 
%pursuing in his own way the works of E. Cartan over the real 
%and of C. Chevalley over an algebraically closed field, has understood, 
%in a highly efficient and non elementary way, the algebraic semisimple groups.

%He has  described  with A. Borel the structure of algebraic semisimple groups 
%over non algebraically closed fields
%in \cite{BorelTits1}, \cite{BorelTits2} and \cite{BorelTits3}.
%He has reduced the classification of  
%these  semisimple groups to the aniso\-tropic casein \cite{Tits66}. 
%He also has 
F. Bruhat and J. Tits have described 
the analog of the Cartan decomposition for  semisimple groups over non-archimedean local fields, in \cite{BruhatTits1},
\cite{BruhatTits2}, \cite{BruhatTits3} and  \cite{BruhatTits4},
by introducing new geometric spaces
that are nowaday called  Bruhat-Tits buildings 
and that extend the space of $p$-adic norms studied 
by Goldman and Iwahori in \cite{GoldmanIwahori}.
\vs 

As explained in the book \cite{KalethaPrasad}, these Bruhat-Tits buildings are very useful.

One of the reason is that they are 
$K(\pi,1)$-spaces for the lattices in semisimple $p$-adic groups.

Another reason is that they played the role of a model to follow
in order to understand other finitely generated groups,
like Coxeter groups, Artin groups, Baumslag-Solitar groups or Mapping class groups.

The relevance of the Bruhat-Tits buildings became even clearer to me when I used them 
with Hee Oh 
to prove a general polar decomposition for $p$-adic symmetric spaces 
in \cite{BenoistOhPolar}. This polar decomposition  was a key ingredient 
in our proof of  equidistribution of $S$-integral points 
on rational symmetric spaces in \cite{BenoistOhSintegral}.

%50
\section{The symplectic adapted basis}
\label{secconcri}

In this section we come back to elementary consideration and we discuss the structure of the {\it rational symplectic group}
${\rm Sp}(n,\m Q)$, and its relation with the  {\it integral symplectic group}
${\rm Sp}(n,\m Z)$. 
\vs

We first recall the well-known undergraduate ``adapted basis theorem'' for $\m Z$-module or, equivalently, the ``Smith normal form'' for integral matrices. 
We denote by ${\mc M}(n,\m Z)$ the ring of $n\times n$ integral matrices.

\bp
\label{prohsidsi0} {\bf (Smith)}
Let $g\in {\mc M}(n,\m Z)$. Then there exist $\si$ and $\si'$ in 
${\rm SL}(n,\m Z)$ and an integral diagonal matrix ${\bf a}={\rm diag}(a_1,\ldots,a_n)$ 
with $a_1|a_2|\ldots|a_{n}$, and  such that 
$$g=\si\;{\bf a}
\;\si'.
$$
\ep

Theorem \ref{thmhsidsi} follows from the  following proposition 
which is a variation of the ``adapted basis theorem'' 
which takes into account the existence of a symplectic form.
We introduce the set ${\mc Mp}(n,\m Z)$ of nonzero integral matrices which are proportional 
to elements of ${\rm Sp}(n,\m R)$,
$$
{\mc Mp}(n,\m Z):=\{ g\in {\mc M}(2n,\m Z)\mid 
{}^tg Jg =\la^2 J
\;\; \mbox{\rm for some $\la$ in $\m R^*$}\}.
$$

\bp
\label{prohsidsi}
Let $g\in {\mc Mp}(n,\m Z)$.
Then there exist two matrices $\si$ and $\si'$ 
in 
${\rm Sp}(n,\m Z)$ and a positive integral diagonal matrix ${\bf a}={\rm diag}(a_1,\ldots,a_{2n})$ 
with  $a_1|a_2|\ldots|a_{n}$, with $a_n|a_{2n}$ and such that
\begin{eqnarray*}
g&=&\si\;
{\bf a}
\;\si'.
\end{eqnarray*}
\ep

Note that the matrix ${\bf a}$ is also in ${\mc Mp}(n,\m Z)$
and hence the products $a_ja_{n+j}$ do not depend 
on the positive integer $j\leq n$. Indeed it is equal to $\la^2$.
In particular, one has $a_{2n}|a_{2n-1}|\ldots |a_{n+1}$.
\vs 

For the proof of Proposition \ref {prohsidsi}, we need the following lemma.
We recall that a nonzero vector $v$  of $\m Z^{k}$ is primitive if it  spans
the $\m Z$-module $\m R v\cap \m Z^k$.

\bl 
\label{lemhsidsi}
The group ${\rm Sp}(n,\m Z)$ acts transitively on the set of primitive vectors in $\m Z^{2n}$.
\el

Denote by $e_1,\ldots, e_n,f_1,\ldots ,f_n$  the  canonical basis of $\m Z^{2n}$ so that our symplectic form is $\om=e_1^*\wedge f_1^*+\cdots+ e_n^*\wedge f_n^*$.

\begin{proof}[Proof of Lemma \ref{lemhsidsi}]
Let $v=(x_1,..,x_{2n})$ be  a
primitive vector in $\m Z^{2n}$.
We want to find $\si\in {\rm Sp}(n,\m Z)$ such that
$\si v=e_1$. 

This is true for $n=1$. Using the subgroups 
${\rm Sp}(1,\m Z)$ for the planes $\m Z e_j\oplus\m Z f_j$, with $j=1,\ldots, n$, we can assume that 
$$
x_{n+1}=\cdots = x_{2n}=0.
$$ 
In this case the vector $(x_1,\ldots,x_n)$
is primitive in $\m Z^n$. 

Since
${\rm SL}(n,\m Z)$ acts transitively on the set of primitive vectors in $\m Z^{n}$, we can find 
a block diagonal matrix 
$\si={\rm diag}(\si_0,{}^t \si_0^{-1})$,
with $\si_0\in {\rm SL}(n,\m Z)$  such that $\si v=e_1$.
This matrix $\si$ belongs to ${\rm Sp}(n,\m Z)$.
\end{proof}

\begin{proof}[Proof of Proposition \ref{prohsidsi}]
Set $\Ga:={\rm Sp}(n,\m Z)$.
The proof is by induction on $n$. It relies on a succession of steps, in the spirit of the Smith normal form,
in which one multiplies on the right or on the left the matrix $g$ by an ``elementary'' matrix to obtain a simpler matrix $g'\in \Ga g\Ga$. We have to pay attention that at each step the elementary matrix is symplectic. 

We can assume that the gcd of the coefficients of $g$ is equal to $1$.
We denote by $\la$ the positive real factor such that $g/\la$ belongs to 
${\rm Sp}(n,\m R)$. Note that $\la^2$ is a positive integer.
At the end of the proof we will see that $a_1=1$ and $a_{n+1}=\la^2$.
\vs 

\noindent
{\bf $\mathbf 1^{\bf st}$ step:} {\it 
We find $g'\in \Ga g\Ga$ such that $g'e_1=e_1$.}
\vs 

Since the coefficients of the integral matrix $ g$ are relatively prime,
by  Proposition \ref{prohsidsi0}, there exists a primitive vector 
$v$ in $\m Z^{2n}$ such that $ gv$ is also primitive.
According to lemma \ref{lemhsidsi}, there exists $\si$, $\si'$ in $\Ga$ such that $\si  g v=e_1$ and  $\si'e_1=v$. 
Then the matrix $g':=\si g\si'$ satisfies $g'e_1=e_1$.
\vs 

\noindent
{\bf $\mathbf 2^{\bf nd}$ step:} 
{\it We find $g'\!\in\! \Ga g\Ga$ with} 
$g'e_1=e_1$ and  $\om(g'e_j,f_1)=0$ {\it for} $j\!>\! 1$.
\vs 

By the first step, we can assume that
\begin{eqnarray*}
g&=&\!\mbox{\scriptsize 
$\left(\!\!\begin{array}{cc} \al\!&\be\\  
\ga\!&\de\end{array}\!\!\right)$}
\;\; \mbox{with $\al e_1=e_1$ and $\ga e_1=  0$}
\end{eqnarray*}
In particular the first column of the integral matrix $\al$ is $(1,0,\ldots,0)$. 
We would like the first row of $\al$ to be also of the form 
$(1,0,\ldots,0)$. For that we choose $g'= g\si'$ where  $\si'$ is the symplectic transformation 
$$\textstyle
\si'=\mathds{1}_n+\sum_{1< j\leq n} \al_{1,j}(f_j\otimes f_1^*-e_1\otimes e_j^*)\in {\rm Sp}(n,\m Z),
$$
where $\al_{1,j}$ are the  coefficients 
of the first row of the matrix $\al$.\vs 

\noindent
{\bf $\mathbf 3^{\bf rd}$ step:} {\it 
We find $g'\in \Ga g\Ga$ such that $g'e_1=e_1$ and $g'f_1=\la^2f_1$.}
\vs 

By the second step, we can assume, writing 
$g=\!\mbox{\scriptsize 
$\left(\!\!\begin{array}{cc} \al\!&\be\\  
\ga\!&\de\end{array}\!\!\right)$}$ that
both the first row and first column of $\al$ are $(1,0,\ldots,0)$, and 
the first column of $\ga$ is $(0,\ldots, 0)$.
We would also like the first row of $\be$ to be 
$(0,\ldots,0)$. For that we choose $g'= g\si'$ where  $\si'$ is the symplectic transformation 
$$\textstyle
\si'=\mathds{1}_n-\be_{1,1}e_1\otimes f_1^*-\sum_{1< j \leq n}\be_{1,j}(e_j\otimes f_1^*+e_1\otimes f_j^*)\in {\rm Sp}(n,\m Z).
$$
Now by construction one  has 
\begin{eqnarray*}
\om(g'e_j,f_1)&=&0
\;\;\;\; {\rm for}\;\;  j<n,\\ 
\om(g'e_1,f_1)&=&\! 1
\;\;{\rm and}\\
\om(g'f_j,f_1)&=&0
\;\;\;\; {\rm for}\;\; j\leq n.
\end{eqnarray*}
Since $g'/\la$ is symplectic, this implies that 
$g'^{-1}f_1=\la^{-2}f_1$, or equivalently, 
$g'f_1= \la^2\, f_1$ as required.
\vs 

\noindent
{\bf $\mathbf 4^{\bf th}$ step:} {\it Conclusion.}
\vs 

By the third step, we can assume that 
$ge_1=e_1$ and $gf_1=\la^{2} f_1$.
Therefore $g$ preserves the symplectic $\m Z$-submodule 
of $\m Z^{2n}$ orthogonal of $\m Ze_1\oplus\m Z f_1$, which admits  $e_2,\ldots,e_{n},f_2,\ldots,f_{n}$ as $\m Z$-basis.
We conclude by applying the induction hypothesis to the restriction 
$g'\in {\mc Mp}(n\!-\!1,\m Z)$ of $g$ to this $\m Z$-module.
\end{proof}

%60
\section{The strong approximation theorem}
\label{secstrapp}

In this section, we give a non elementary proof of 
the decomposition theorem \ref{thmhsidsi} for 
${\rm Sp}(n,\m Q)$.  We will deduce this theorem from
the Bruhat-Tits decomposition theorem \ref{thmhsidsik} for 
${\rm Sp}(n,\m Q_p)$ thanks to the strong approximation theorem. 
\vs 

First, I recall the strong approximation theorem.
I will not quote here 
the general formulation for a simply-connected isotropic $\m Q$-simple 
algebraic group defined over $\m Q$ that can be found in \cite{PlatonovRapinchuk94}
but only the special case where $G$ is the symplectic group.
\vs 

For $p=2,3,5,\ldots$ a prime number, we denote by 
$\m Q_p$ the $p$-adic local field and by $\m Z_p$ 
its ring of integers. 
%which is its maximal open compact subring.

We denote by $\widehat{\m Q}=\prod'_p\m Q_p$ the locally compact ring of finite
ad\`eles which is the restricted product of the $\m Q_p$ with respect to the open subrings  $\m Z_p$. The product
$\widehat{\m Z}:=\prod_p \m Z_p$ is then a maximal compact open subring of $\widehat{\m Q}$.   

Note that, thanks to the diagonal embedding,  
$\m Q$ is a dense subring in $\widehat{\m Q}$.
This means that $\widehat{\m Q}=\m Q +\widehat{\m Z}$ and that 
$\m Z$ is dense in $\widehat{\m Z}$.
%This means that, for every family $x=(x_p)$ with $x_p\in \m Q_p$,
%such that, for almost all $p$, $x_p$ is an integer,
%and, for all integers $N\geq 1$, 
%there exists $y\in \m Q$ such that 
%$x\!-\!y$ belongs to the compact open subring $N\widehat{\m Z}$.
\vs 

By construction the symplectic group 
${\rm Sp}(n,\widehat{\m Q})$ is a locally compact group 
that contains ${\rm Sp}(n,\widehat{\m Z})$ as a maximal compact open
subgroup. It also contains the group ${\rm Sp}(n,\m Q)$.

Here is the strong approximation theorem for the symplectic group.

\bt
\label{thmstrapp} 
The group ${\rm Sp}(n,\m Q)$ is dense in ${\rm Sp}(n,\widehat{\m Q})$.
\et

This means that, 
$${\rm Sp}(n,\widehat{\m Q})= {\rm Sp}(n,\m Q){\rm Sp}(n,\widehat{\m Z})$$ 
and that 

\centerline{${\rm Sp}(n,\m Z)$ is dense in ${\rm Sp}(n,\widehat{\m Z})$.}
\vs

%for every family $g=(g_p)$ with $g_p\in {\rm Sp}(n,\m Q_p)$,
%such that, for almost all $p$, $g_p$ has integer coefficients,
%and, for all integers $N\geq 1$, 
%there exists $h\in {\rm Sp}(n,\m Q)$ such that 
%$gh^{-1}$ belongs to the compact open subgroup $K_N$
%defined by
%$$K_N:=\{k\in {\rm Sp}(n,\widehat{\m Z})\mid k\equiv{\bf 1}\; {\rm mod} N\}$$

If we collect together the Bruhat-Tits decomposition
in Theorem \ref{thmhsidsik} for all $p$-adic fields $k=\m Q_p$,
one gets

\bt
\label{thmhsidsia} 
Let $g\in {\rm Sp}(n,\widehat{\m Q})$. Then there exist two matrices $\si$ and $\si'$ in 
${\rm Sp}(n,\widehat{\m Z})$ and a positive integral diagonal matrix ${\bf d}={\rm diag}(d_1,\ldots,d_n)$ 
with $d_1|d_2|\ldots |d_n$  such that
\begin{eqnarray*}
	g&=&\si\;
	\mbox{\small%\scriptsize 
		$\left(\!
		\begin{array}{cc} 
			{\bf d}&{\bf 0}\\  
			{\bf 0}&\!{\bf d}^{-1}\!
		\end{array}
		\!\right)$}
	\;\si'.
\end{eqnarray*}
\et

We can now give the non-elementary proof of the symplectic 
Smith normal form.

\begin{proof}[Proof of Theorem \ref{thmhsidsi}]
Let $g\in {\rm Sp}(n,\m Q)$. 

According to the combined Bruhat-Tits decomposition  theorem \ref{thmhsidsia},
one can write 
$$
g=\si\,{\bf a}\,\si'
$$ 
with
$\si$, $\si'$ in 
${\rm Sp}(n,\widehat{\m Z})$ 
and with 
${\bf a}=\mbox{\small%\scriptsize 
	$\left(\!
	\begin{array}{cc} 
		{\bf d}&{\bf 0}\\  
		{\bf 0}&\!{\bf d}^{-1}\!
	\end{array}
	\!\right)$}$
where 
${\bf d}={\rm diag}(d_1,\ldots,d_n)$ is a positive integral diagonal matrix 
with $d_1|d_2|\ldots |d_n$.

According to the strong approximation theorem \ref{thmstrapp}, 
one can write 
$$\si=\si_0\eta$$ 
with $\si_0$ in ${\rm Sp}(n,\m Z)$ and with 
$\eta$ in an arbitrarily small neighborhood of ${\bf 1}$ 
in  ${\rm Sp}(n,\widehat{\m Z})$. 
More precisely we choose $\eta$ such that the element 
$\si'_0:={\bf a}^{-1}\eta\,{\bf a}\,\si'$ belongs to 
 ${\rm Sp}(n,\widehat{\m Z})$.
Then one has the equality
$$
g=\si_0\,{\bf a}\,\si'_0
$$ 
where both
$\si_0$ and $\si'_0={\bf a}^{-1}\si_0^{-1}g$ belong to
${\rm Sp}(n,\m Z)$.
\end{proof}

%70
\section{Chevalley groups}
\label{secchegro} 

Let ${\bf G}$ be a simply-connected 
Chevalley group.  See 
\cite{Steinberg68} for a concrete presentation of the group ${\bf G}(\m Z)$,
and see \cite{Gross96} for other nice examples of
$\m Z$-models of simple algebraic groups over $\m Q$.
This ${\bf G}$ is a reductive scheme-group over $\m Z$ 
such that as a $\m Q$-group ${\bf G}$ is a $\m Q$-split simply connected quasi-simple algebraic group.
By construction, this algebraic group contains a $\m Q$-split maximal torus 
${\bf A}$ such that the group of integral points ${\bf N}(\m Z)$ of the normalizer 
of ${\bf A}$ surjects onto the Weyl group
of $({\bf A}(\m C),{\bf G}(\m C)$.

Since ${\bf G}$ is simply connected and $\m R$-isotropic, by strong 
approximation,
the group ${\bf G}(\m Q)$ is dense in 
${\bf G}(\widehat{\m Q})$.
On the other hand, for all prime integer $p$,
one can consider 
the simply connected simple 
$p$-adic Lie group $G:={\bf G}(\m Q_p)$, 
its split maximal torus $A:={\bf A}(\m Q_p)$
and its normalizer $N:={\bf N}(\m Q_p)$.
The maximal compact subgroup $K:={\bf G}(\m Z_p)$ is a good compact subgroup in the sense that one has the equality
$N=(N\cap K)A$. 
Hence, according to Bruhat-Tits, one has the decomposition
${\bf G}(\m Q_p)={\bf G}(\m Z_p){\bf A}(\m Q_p){\bf G}(\m Z_p)$.

Therefore the same proof as in 
Chapter \ref{secstrapp} gives the following theorem due to Steinberg in
\cite[Theorem 21]{Steinberg68}

\bt
\label{thmhsidsic}
Let ${\bf G}$ be a simply connected Chevalley group and $g\in {\bf G}(\m Q)$.
Then there exist two elements
$\si$ and $\si'$ 
in 
${\bf G}(\m Z)$ and an element ${\bf a}$
in ${\bf A}(\m Q)$ such that
\begin{eqnarray*}
g&=&\si\;
{\bf a}
\;\si'.
\end{eqnarray*}
\et

\noindent
{\bf Remark}. 
Such a decomposition is not true 
when we replace $\m Q$ by a number field
$\m K$ whose ring of integer $\mc O$
is not principal. 
Here is an example with
$${\bf G}(\m K):={\rm SL}(2,\m K)
\; ,\;\;
{\bf G}(\mc O):={\rm SL}(2,\mc O)\; ,$$ 
\begin{eqnarray*}
{\bf A}(\m K)&:=&\left\{ {\bf a}=\;
\mbox{\small%\scriptsize 
$\left(\!
\begin{array}{cc} 
{\bf d}&{\bf 0}\\  
{\bf 0}&\!{\bf d}^{-1}\!
\end{array}
\!\right)$}
\;\mid {\bf d}\in K^*\right\}.
\end{eqnarray*}
In this case the product
${\bf G}(\mc O)\,{\bf A}(\m K)\,{\bf G}(\mc O)$
is not equal to ${\bf G}(\m K)$. 
For instance, when $K=\m Q[i\sqrt{5}]$
and $\mc O=\m Z[i\sqrt{5}]$, this product does not contain the matrix
\begin{eqnarray*}
g=\;
\mbox{\small%\scriptsize 
$\left(\!
\begin{array}{cc} 
\!(1\!-\! i\sqrt{5})/2\!&i\sqrt{5}\\  
-1&\! 2\!
\end{array}
\!\right)$}.
\end{eqnarray*}
Indeed the element ${\bf d}\in K^*$ should be a unit in all completion $K_{\g p}$ except for the prime ideal $\g p_0=2\m Z\oplus (1\!+\!i\sqrt{5})\m Z$ in which case 
it should be a uniformizer. 
Such an element ${\bf d}$ would be a generator
of the ideal $\g p_0$. This is a contradiction, since this ideal $\g p_0$ is not principal.

%9
{\small
\bibliography{symplecticgroup}
}
\vs 

{\small
\noindent
Y. \textsc{Benoist}: CNRS, 
Universit\'e Paris-Saclay,\hfill
%e-mail: 
\texttt{yves.benoist@u-psud.fr}}

\end{document}